\documentclass[11pt,3p,authoryear]{elsarticle}
\usepackage{natbib,graphics,amsmath,amsfonts,amssymb,tabularx,rotating,geometry,multirow,diagbox,xcolor, longtable,enumerate,amsthm,hyperref,listings}
\usepackage{hyphenat}
\usepackage[english]{babel}
\newcommand{\pcite}[1]{\citeauthor{#1}'s (\citeyear{#1})}

\journal{arXiv}

\title{{Nonfractional Memory:\\ Filtering, Antipersistence, and Forecasting}}

\author{J. Eduardo Vera-Vald\'es\corref{fn1}}
\address{Department of Mathematical Sciences, Aalborg University, and CREATES.}
\cortext[fn1]{E-mail: \href{mailto:eduardo@cimat.mx}{eduardo@math.aau.dk} Webpage: \url{https://sites.google.com/view/veravaldes/} \\ Address: Skjernvej 4A, 9220 Aalborg East, Denmark. }

\newtheorem{teo}{Theorem}

\allowdisplaybreaks

\begin{document}

\begin{abstract}
\noindent The fractional difference operator remains to be the most popular mechanism to generate long memory due to the existence of efficient algorithms for their simulation and forecasting. Nonetheless, there is no theoretical argument linking the fractional difference operator with the presence of long memory in real data. In this regard, one of the most predominant theoretical explanations for the presence of long memory is cross-sectional aggregation of persistent micro units. Yet, the type of processes obtained by cross-sectional aggregation differs from the one due to fractional differencing. Thus, this paper develops fast algorithms to generate and forecast long memory by cross-sectional aggregation. Moreover, it is shown that the antipersistent phenomenon that arises for negative degrees of memory in the fractional difference literature is not present for cross-sectionally aggregated processes. Pointedly, while the autocorrelations for the fractional difference operator are negative for negative degrees of memory by construction, this restriction does not apply to the cross-sectional aggregated scheme. We show that this has implications for long memory tests in the frequency domain, which will be misspecified for cross-sectionally aggregated processes with negative degrees of memory. Finally, we assess the forecast performance of high-order $AR$ and $ARFIMA$ models when the long memory series are generated by cross-sectional aggregation. Our results are of interest to practitioners developing forecasts of long memory variables like inflation, volatility, and climate data, where aggregation may be the source of long memory.

\bigskip

\noindent {\it Keywords:} Nonfractional memory, long memory, fractional difference, antipersistence, forecasts.

\bigskip
\noindent {\it JEL classification:} C15, C22, C53.

\end{abstract}

\maketitle

\clearpage

\section{Introduction}

Long memory has been a topic of interest in econometrics since \pcite{Granger1966} study on the shape of the spectrum of economic variables. Granger found that {\it long-term fluctuations in economic variables if decomposed into frequency components are such that the amplitudes of the components decrease smoothly with decreasing period}. As shown by \cite{Adenstedt1974}, this type of behaviour implies long lasting autocorrelations. The presence of long memory in the data can have perverse effects in estimation and forecasting methods if not included into the modelling scheme, see \cite{Beran1994}, and \cite{Beran2013}.

The autoregressive fractionally integrated moving average $(ARFIMA)$ class of models has become one of the most popular methods to model long memory in the time series literature. They have the appeal of bridging the gap between the stationary autoregressive moving average $(ARMA)$ models and the nonstationary autoregressive integrated moving average $(ARIMA)$ model. $ARFIMA$ models rely on the fractional difference operator to introduce long memory behaviour. Nonetheless, there are currently no economic or financial reasonings implying the fractional difference operator. 

One of the reasons behind the reliance of the time series literature on generating long memory by the fractional difference operator is the existence of efficient algorithms for their simulation and forecasting. In general, these type of algorithms are not available for other long memory generating schemes. Thus, in this paper, we develop algorithms for memory generation and forecasting by cross-sectional aggregation \`a la \cite{Granger1980}. 

We then contrast the properties of cross-sectionally aggregated processes to the ones obtained by the fractional difference operator. It is shown that cross-sectional aggregated processes are more flexible than fractionally differenced ones, in the sense that more short term dynamics may be included. 

Moreover, we show that the antipersistent property of negative autocorrelation for negative degrees of memory does not apply to long memory generated by cross-sectional aggregation. We show that this has repercussions for long memory estimators in the frequency domain. 

Finally, this paper evaluates the forecasting power of the $ARFIMA$, and high-order $AR$ models when forecasting long memory generated by cross-sectional aggregation. It finds that high-order $AR$ models beat a pure fractional difference process, $I(d)$, in terms of forecasting performance for some cases. Nonetheless, allowing for short term dynamic in the form of an $ARFIMA(1,d,0)$ model produces comparable forecast performance as high-order $AR$ models while relying on fewer parameters. 

This paper proceeds as follows. In Section \ref{Sec:ARFIMA}, we present the fractional difference operator commonly used to model long memory in the time series literature. Section \ref{Sec:CSA} discusses cross-sectional aggregation as the theoretical explanation behind the presence of long memory in the data and develops a fast algorithm for its generation. Moreover, it contrasts properties of the fractional difference operator against the cross-sectional aggregation scheme. Section \ref{Sec:Ant} discusses the antipersistent property in the context of cross-sectionally aggregated processes. Section \ref{Sec:For} constructs minimum squared error forecasts, and studies the theoretical performance of high-order $AR$, and $ARFIMA$ models when forecasting long memory generated by cross-sectional aggregation. Section \ref{Sec:Con} concludes. 

\section{The Fractional Difference Operator}\label{Sec:ARFIMA}

The $ARFIMA$ specification due to \cite{Granger1980b}, and \cite{Hosking1981} has become the standard model to study long memory in the time series literature. They extended the $ARMA$ model to include long memory dynamics by introducing the fractional difference operator
\begin{equation*}
(1-L)^d x_t = \varepsilon_t,
\end{equation*}
where $\varepsilon_t$ is a white noise process, and $d\in(-1/2,1/2)$. Following the standard binomial expansion, they decompose the fractional difference operator, $(1-L)^d$, to generate a series given by
\begin{equation}\label{eq:arfima_ma}
x_t = \sum_{j=0}^{\infty} \pi_j \varepsilon_{t-j},
\end{equation} 
with coefficients $\pi_j=\Gamma(j+d)/(\Gamma(d)\Gamma(j+1))$ for $j\in\mathbb{N}$. Using Stirling's approximation, it can be shown that these coefficients decay at a hyperbolic rate, $\pi_j\approx j^{d-1}$, which in turn translates to slowly decaying autocorrelations. We write $x_t\sim I(d)$ to denote a process generated by equation (\ref{eq:arfima_ma}), that is, an integrated process with parameter $d$.

The properties of the $ARFIMA$ model have been well documented in, among others, \cite{Baillie1996}, and \cite{Beran2013}. Moreover, fast algorithms have been developed to generate series using the fractional difference operator, see \cite{Jensen2014}. Thus, the $ARFIMA$ model has become the canonical construction for modelling and forecasting long memory in the time series literature. 

Even though the $ARFIMA$ provides a good representation of long memory, and bridges the gap between the stationary $ARMA$ models and the non-stationary $ARIMA$ model, to the best of our knowledge, there is no theoretical reasoning linking the $ARFIMA$ model with the memory found in real data. That is, in contrast with the complete market hypothesis implying that stock prices follow a random walk, or capital depreciation suggesting an autoregressive process, there are no economic or financial arguments for the fractional difference operator to occur in real data. In this regard, the next section presents the most common theoretical motivation behind long memory in the time series literature, cross-sectional aggregation.

\section{Cross-Sectional Aggregation}\label{Sec:CSA}

\cite{Granger1980}, in line with the work of \cite{Robinson1978} on autoregressive processes with random coefficients, showed that aggregating $AR(1)$ processes with coefficients sampled from a Beta distribution can produce long memory. 

Define a cross-sectional aggregated series as
\begin{equation}\label{eq:csa_def}
x_t = \frac{1}{\sqrt{N}}\sum_{i=1}^N x_{i,t},
\end{equation}
where the $N$ individual series are generated as
$$x_{i,t} = \alpha_i x_{i,t-1}+\epsilon_{i,t}\ \ \ i=1, 2, \cdots, N;$$
where $\varepsilon_{i,t}$ is an $i.i.d.$ process with $E[\epsilon_{i,t}^2] = \sigma_\epsilon^2$. Moreover, $\alpha_i^2$ is sampled from a Beta distribution with parameters $a,b>0$, with density given by
$$\mathcal{B}(\alpha; a, b) =  \frac{1}{B(a,b)} \alpha^{a-1}(1-\alpha)^{b-1}\ \ \ \ \text{for}\ \ \ \alpha\in(0,1),$$
where $B(\cdot,\cdot)$ is the Beta function. 

Granger showed that, as $N\to\infty$, the autocorrelation function of $x_t$ decays at a hyperbolic rate with parameter $d=1-b/2$. Thus, $x_t$ has long memory.

The cross-sectional aggregation result has been extended in several directions, including to allow for general $ARMA$ processes, as well as to other distributions. See for instance, \cite{Linden1999}, \cite{Oppenheim2004}, and \cite{Zaffaroni2004}. As argued by \cite{Haldrup2017}, maintaining the Beta distribution allows us to have closed form representations. Furthermore, \cite{Beran2010} proposed a method to estimate the parameters from the Beta distribution from the individual series, $x_{i,t}$ above.

In applied work, cross-sectional aggregation has been cited as source of long memory for many series. To name but a few, it has been proposed for inflation, output, and volatility; see \cite{Balcilar2004}, \cite{Diebold1989}, and \cite{Osterrieder2015}.  

Even though a cross-sectional aggregated process has long memory, \cite{Haldrup2017} show that the resulting series does not belong to the $ARFIMA$ class of processes. Pointedly, they show that the fractionally differenced cross-sectional aggregated process does not follow an $ARMA$ specification. 

The next subsections expands on the cross-sectional aggregation literature by developing fast algorithms for its generation, and contrasting its theoretical properties to the ones of the fractional difference operator.

\subsection{Nonfractional Memory Generation}

As can be seen by its definition, Equation (\ref{eq:csa_def}), generating cross-sectionally aggregated processes is computationally demanding. For each cross-sectionally aggregated process, we need to simulate a vast number of $AR(1)$ processes, which are in turn computationally demanding. \cite{Haldrup2017} suggest that to get a good approximation, the cross-sectional dimension should increase as the sample size. The computational demands are thus particularly large for Monte Carlo type of analysis on cross-sectionally aggregated processes. In what follows, we use the theoretical autocorrelations of the cross-sectionally aggregated process to present an algorithm that makes generation of long memory by cross-sectional aggregation comparable to fractional differencing in terms of computational requirements. 

Denote $x_t\sim CSA(a,b)$ to a series generated by cross-sectional aggregation of autoregressive parameters sampled from the Beta distribution $B(a,b)$, Equation (\ref{eq:csa_def}). The notation makes explicit the origin of the memory by cross-sectional aggregation and its dependence on the two parameters of the Beta distribution. Theorem \ref{teo:csa_gen} develops a fast algorithm to generate cross-sectionally aggregated processes.

\begin{teo}\label{teo:csa_gen}
	Let $x_t\sim CSA(a,b)$ defined as in (\ref{eq:csa_def}), then $x_t$ can be computed as the first $T$ elements of the $(2T-1)\times 1$ vector 
	$$T^{-1} F(\bar{F}\tilde{\phi}\odot \bar{F}\tilde{\varepsilon}),$$
	where $\bar{F}$ is the discrete Fourier transform, $T^{-1}F$ is the inverse transform, and '$\odot$' denotes multiplication element by element. Furthermore, $\tilde{z}$ is a $(2T-1)\times 1$ vector given by $\tilde{z}:=[z_0, z_1, \cdots, z_{T-1},$ $0, \cdots, 0]$ for $z_t=\phi_t,\varepsilon_t$ where $\varepsilon_j\sim i.i.d. N(0,\sigma^2_\epsilon)$, and $\phi_j = \left(B(p+j,q)/B(p,q)\right)^{1/2},\ \forall j\in\mathbb{N}$.
\end{teo}
\noindent Proof: See appendix \ref{app:proofs}.

The theorem is an application of the circular convolution theorem for the coefficients associated with the cross-sectionally aggregated process. In this sense, it is in line with the algorithm of \cite{Jensen2014} for the fractional difference operator and thus achieves equivalent computational efficiency. The difference in computation times between the standard simulation of a cross-sectional aggregated times series, Equation (\ref{eq:csa_def}), and the fast implementation in Theorem \ref{teo:csa_gen} are very large for all sizes. For small series $(T\approx 10^2)$ the gains are in the hundreds of times faster, while in the thousands for medium sized series $(T\approx 10^4)$. This is of course not surprising given that, for a sample of size $T$, the number of computations needed for the standard implementation is of order $NT^{2}$, with $N$ the cross-sectional dimension, which, as argued before, should increase as $T$ does. Meanwhile, the computational requirements for the fast implementation are of order $T \log{(T)}$. Codes implementing the algorithm for memory generation by cross-sectional aggregation are available in Appendix \ref{app:codes}, and on the author's github repository at \href{https://github.com/everval}{github.com/everval/Nonfractional}.

To get a better understanding of the dynamics of the long memory by cross-sectional aggregation, the next sections compare it against the pure fractional noise, that is, a fractionally differenced white noise. However, note that we can allow more short-term dynamic by adding $AR$ and $MA$ filters to both specifications.

\subsection{Nonfractional Memory and Fractional Difference Operator}

One notable difference between a cross-sectionally aggregated process and a fractionally differenced one is the number of parameters needed for its generation. Fractionally differenced processes rely on only one parameter to model the entire series, whilst a cross-sectionally aggregated process uses two. As we will see below, this gives the cross-sectional aggregation procedure more flexibility. 

Let $\gamma_{I(d)}(\cdot)$, and $\gamma_{CSA(a,b)}(\cdot)$ be the autocorrelation functions of an $I(d)$, and a $CSA(a,b)$ process, respectively. That is,
\begin{equation}\label{ACF:Id}
\gamma_{I(d)}(k) = \frac{\Gamma(k+d)\Gamma(1-d)}{\Gamma(k-d+1)\Gamma(d)} \approx k^{2d-1}\ \ \ \ \text{ as } k\to\infty
\end{equation}
\begin{equation}\label{ACF:CSA}
\gamma_{CSA(a,b)}(k) = \frac{B(a+k/2,b-1)}{B(a,b-1)} \approx k^{1-b}\ \ \ \text{as } k\to\infty,
\end{equation}
which show that both processes have hyperbolic decaying autocorrelations. As argued by \cite{Granger1980}, the asymptotic behaviour of the autocorrelation function for a cross-sectionally aggregated process, Equation (\ref{ACF:CSA}), only depends on the second argument of the Beta function. In this context, both processes show the same rate of decay for the autocorrelation function for $b=2(1-d)$. 

\begin{figure}[h!]
	\centering
	\includegraphics[scale=1]{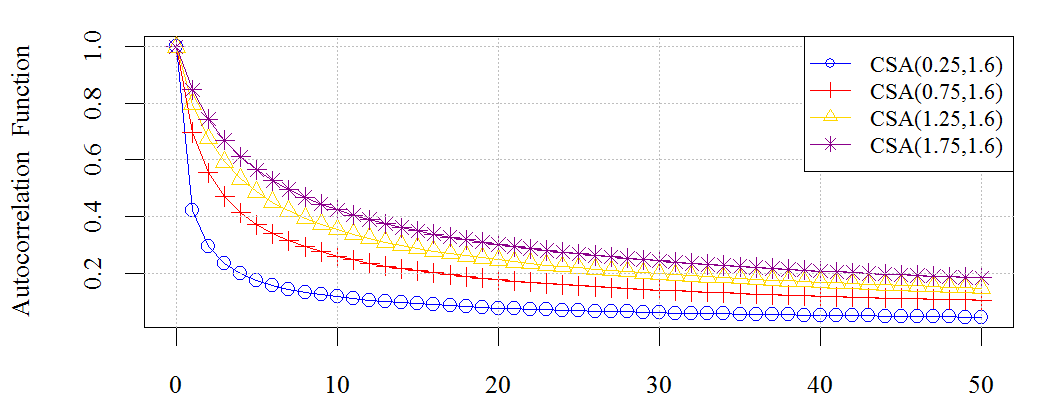}  
	\caption{Autocorrelation function for an $CSA(a,b)$ processes for different values of `$a$' while having the same asymptotic behaviour.}
	\label{fig:ACF_CSA_ps}
\end{figure}

Moreover, the first argument of the Beta function allows us to introduce short term dynamics. Figure \ref{fig:ACF_CSA_ps} shows the autocorrelation function of cross-sectional aggregated processes for different values of `$a$', the first parameter of the Beta distribution. The figure shows that as the first parameter gets larger, so does the autocorrelation function for the initial lags. Thus, `$a$' acts as a short memory modulator.

Section \ref{Sec:For} will look at the forecasting performance of $ARFIMA$ models when working with cross-sectionally aggregated processes, but for now suppose we are interested in looking in the other direction. That is, assume we have a process generated by the fractional difference operator and we want to approximate it with a cross-sectionally aggregated process. We show that we can use the first argument to generate such a series. 

Consider the loss function 
\begin{equation}\label{loss:csa2fi}
\mathcal{L}(k,a,d) := \sum_{i=0}^{k}\left(\gamma_{I(d)}(i)-\gamma_{CSA(a,2(1-d))}(i)\right)^2,
\end{equation}
which measures the squared difference between autocorrelations at the first $k$ lags for cross-sectionally aggregated and fractionally differenced processes with the same long memory dynamics, hence we set $b=2(1-d)$ for the cross-sectionally aggregated process.

Minimizing $(\ref{loss:csa2fi})$ with respect to `$a$', allows us to find the cross-sectional aggregated process that best approximates an $I(d)$ one up to lag $k$, while having the same long term dynamics. Given the different forms of the autocorrelation functions, (\ref{ACF:Id}), and (\ref{ACF:CSA}), there is in general not a value of `$a$' that minimizes $(\ref{loss:csa2fi})$ for all values of $k$; for instance, $\min_{a}\mathcal{L}(2,a,0.2)=0.118$, while $\min_{a}\mathcal{L}(30,a,0.2)=0.121$. Moreover, it will depend on $d$.  

\begin{figure}[h!]
	\centering
	\includegraphics[scale=1]{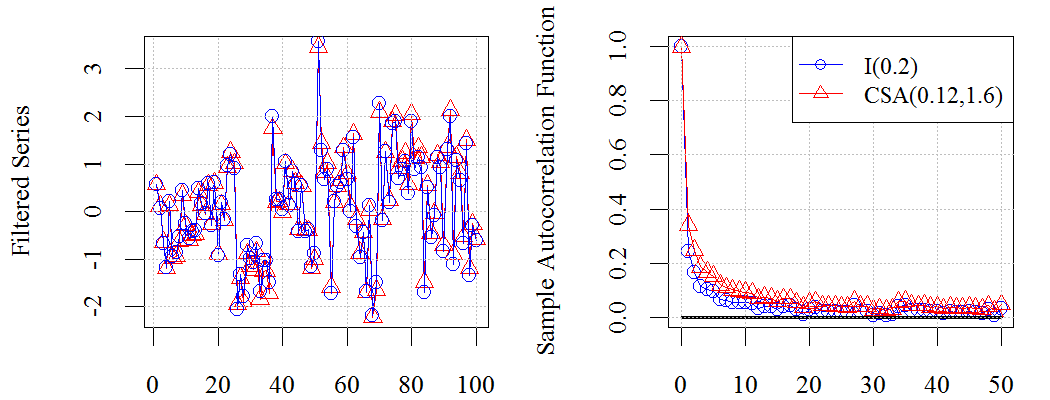}  
	\caption{Filtered series and autocorrelation function for a fractional differenced process $I(0.2)$, and a cross-sectional aggregated one $CSA(0.12,1.6)$.}
	\label{fig:ACF_CSA_Filter}
\end{figure}

Nonetheless, selecting a medium sized $k$, say $k\approx 10$, the approximation turns out to be quite good in general. In Figure \ref{fig:ACF_CSA_Filter}, we present a filtered $\{\varepsilon_t\}_{t=1}^{10^4}\sim N(0,1)$ vector using the fractional difference operator with long memory parameter $d=0.2$, and using the cross-sectional aggregated algorithm with parameters $a=0.12$, and $b=1.6$, so that they show similar short and long memory behaviours. 

The figure shows that the filtered series are quite similar. This behaviour can further be seen on the autocorrelation functions showing similar dynamics. Thus, the figure shows that it is possible to generate cross-sectionally aggregated processes that closely mimic ones due to fractional differencing.

\section{Nonfractional Memory and Antipersistent Processes}\label{Sec:Ant}

It is well known in the long memory literature that the fractional difference operator implies that the autocorrelation function is negative for negative degrees of memory, $d\in (-1/2,0)$. This can be seen in (\ref{ACF:Id}) where the sign of $\gamma_{I(d)}(k)$ depends on $\Gamma(d)$ in the denominator, which is negative for $d\in (-1/2,0)$. Furthermore, the behaviour of the spectral density for a fractionally differenced process near the origin is given by
\begin{equation}\label{spec:id}
f_{I(d)}(\lambda) \approx c_0 \lambda^{-2d}\ \ \ \text{as} \ \ \ \lambda\to 0,
\end{equation}
where $c_0$ is a constant. Thus, $f_{I(d)}(\lambda) \to 0$ as $\lambda\to 0$, that is, the fractional difference operator for negative degrees of memory imply a spectral density collapsing to zero at the origin. These properties, among other related ones, has been named {\it antipersistence} in the literature.

\begin{figure}[h!]
	\centering
	\includegraphics[scale=1]{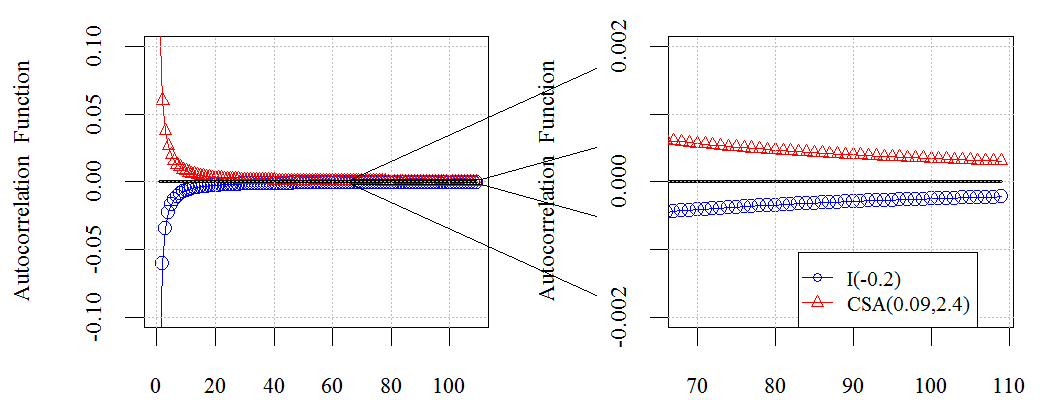} 
	\caption{Autocorrelation functions for a fractional differenced process, $I(-0.2)$, and a cross-sectional aggregated one, $CSA(0.09,2,4)$. The right plot shows a close-up for lags 70 to 110.}
	\label{fig:ACF_CSA_2_FI}
\end{figure}

We argue that cross-sectionally aggregated processes do not share these features. First, Equation (\ref{ACF:CSA}) shows that the autocorrelation function for the cross-sectionally aggregated process only depends on the Beta function, which is always positive. Figure \ref{fig:ACF_CSA_2_FI} shows the autocorrelation function for both fractionally differenced and cross-sectional aggregated processes for a negative degree of memory, $d=-0.2$. The figure shows that, even though both processes show the same rate of decay in their autocorrelation functions, they have opposite signs. 

Then, given the positive sign for all autocorrelations, Theorem \ref{teo:ant} shows that the spectral density for the cross-sectionally aggregated process for negative degrees of memory converges to a positive constant. 

\begin{teo}\label{teo:ant}
	Let $x_t\sim CSA(a,b)$ defined as in (\ref{eq:csa_def}) with $b\in(2,3)$ so that the long memory parameter is in the negative range, $d\in(-1/2,0)$, then, the spectral density of $x_t$ at the origin is positive. That is,
	\begin{equation}\label{spec:csa}
	f_{CSA(a,b)}(0)=c_{a,b}>0,
	\end{equation}
	where $c_{a,b}$ depends on the parameters of the Beta distribution.
\end{teo}
\noindent Proof: See appendix \ref{app:proofs}.

Figure \ref{fig:Ant_Spec_ds} shows the periodogram, an estimate of the spectral density, for $CSA(a,b)$ and $I(d)$ processes of size $T=10^4$ averaged for $10^4$ replications. The figure shows that the periodogram for both processes show similar patterns for positive degrees of memory, both diverging to infinity at the same rate. Nonetheless, for negative degrees of memory, the periodogram collapses to zero as the frequency goes to zero for the fractionally differenced process, while it converges to a constant for the cross-sectionally aggregated process. 

\begin{figure}[h!]
	\centering
	\includegraphics[scale=1]{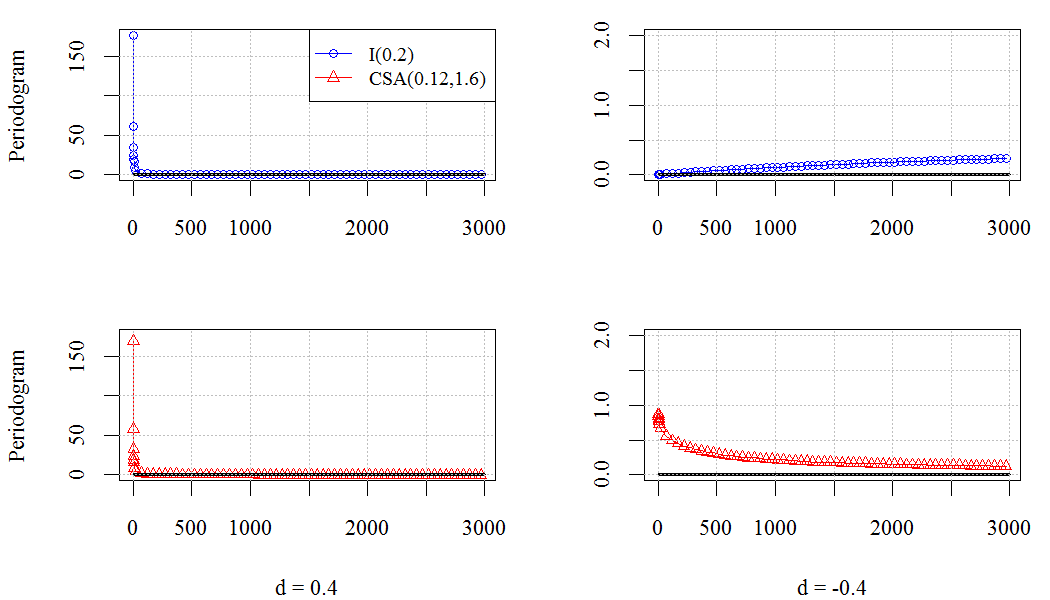} 
	\caption{Mean periodogram of a fractional differenced process, $I(d)$, [top], and a cross-sectional aggregated process, $CSA(0.2,2(1-d))$, [bottom]. A sample size of $T=10,000$ was used, with 10,000 replications, and $d=0.4,-0.4$.}
	\label{fig:Ant_Spec_ds}
\end{figure}

Moreover, the latter property has implications for estimation and inference. In particular, tests for long memory in the frequency domain will be affected. These tests are based on the rate to which the periodogram goes to zero as an estimator for long memory using the log-periodogram regression, see \cite{Geweke1983}, and \cite{Robinson1995a}. 

The log-periodogram regression is given by 
$$\log(\hat{f}_X(\lambda_j)) = a-2d \log(\lambda_j)+u_j,$$
where $\hat{f}_X(\cdot)$ is the periodogram, $a$ is a constant, and $u_j$ is the error term. From (\ref{spec:id}), note that the log-periodogram regression provides an estimate of the long memory parameter, $d$, for fractionally differenced processes. Tests in the frequency domain use this expression to estimate the degree of memory. Nonetheless, as Theorem \ref{teo:ant}, and Figure \ref{fig:Ant_Spec_ds} show, these tests will be misspecified for long memory by cross-sectional aggregation for negative degrees of memory.

To illustrate the misspecification problem, Table \ref{tab:Ant_Spec} reports the degree of long memory estimated by the method of \citealp{Geweke1983}, $GPH$ hereinafter, for several degrees of memory for both fractionally differenced, and cross-sectionally aggregated processes. 

\begin{table}[ht!]
	\begin{small}
		\setlength{\tabcolsep}{3pt}
		\renewcommand{\arraystretch}{1.2} 
		\caption{Mean, and standard deviation in parentheses, of estimated long memory by the $GPH$ method for $CSA(0.2,2(1-d))$ and $I(d)$ processes. We use the standard $T^{1/2}$ bandwidth, where the sample size is $T=10^4$, and with $10^4$ replications.}
		\begin{center}
			\begin{tabular}{cc|cc|cc|cc}
				\multicolumn{2}{c|}{$d=0.4$}&\multicolumn{2}{c|}{$d=0.2$}&\multicolumn{2}{c|}{$d=-0.2$}&\multicolumn{2}{c}{$d=-0.4$}\\
				\hline
				$CSA(a,b)$ &$I(d)$ &$CSA(a,b)$ &$I(d)$ &$CSA(a,b)$ &$I(d)$ &$CSA(a,b)$ &$I(d)$ \\
				0.4062   & 0.4034  &  0.2628   & 0.2011  &  0.1036  & -0.1985  &  0.0653  & -0.3927\\
				(0.0701)   &( 0.0697)  &(0.0697) &(0.0696) &(0.0687) &(0.0685)  &(0.0695)  &(0.0717)\\
				\hline
			\end{tabular}\label{tab:Ant_Spec}
		\end{center}
	\end{small}
\end{table}

The table shows that the estimator is relatively close to the true memory for both processes when the memory is positive, if slightly overshooting it for the cross-sectionally aggregated process, as reported by \cite{Haldrup2017}. This contrasts to the case of negative memory, where the table shows that the estimator remains precise for the fractionally differenced series, while it is incapable of detecting the long memory in the cross-sectionally aggregated processes, with estimators not statistically different from zero. This is of course not surprising in light of Theorem \ref{teo:ant}. 

In sum, the lack of the antipersistent property in cross-sectionally aggregated processes shows that care must be taken when estimating long memory for negative degrees of memory if the long memory generating mechanism is not the fractional difference operator. Further analysis of negative degrees of memory and antipersistence is a line of inquiry open for further research.

\section{Nonfractional Memory Forecasting}\label{Sec:For}

In the time series literature, some effort has been directed to assess the performance of the $ARFIMA$ type of models when forecasting long memory processes. For instance, \cite{Ray1993} calculates the percentage increase in mean-squared error ($MSE$) from forecasting $I(d)$ series with $AR$ models. She argues that the $MSE$ may not increase significantly, particularly when we do not know the true long memory parameter. \cite{Crato1996} compare the forecasting performance of $ARFIMA$ models against $ARMA$ alternatives and find that $ARFIMA$ models are in general outperformed by $ARMA$ alternatives. Moreover, \cite{Man2003} argues that an $ARMA(2,2)$ model compares favourably to an $I(d)$ for short-term forecasts of long memory time series with fractionally differenced structure.

One thing that these forecasting comparison studies have in common is the underlying assumption that long memory is generated by an $ARFIMA$ process. In this context, \textit{a priori} the forecasting exercises assume that a fractionally differenced process is the correct specification. As previously discussed, even though cross-sectional aggregation does generate long memory, the series does not follow an $ARFIMA$ specification. The question remains whether an $ARFIMA$ model serves as a good approximation for forecasting purposes.

The next subsection computes the minimum square error forecasts for cross-sectional aggregated processes, $CSA(a,b)$. With those as benchmark, the following subsections evaluate the forecasts performance of high order $AR$, and $ARFIMA$ models on $CSA(a,b)$ processes.

\subsection{Minimum Square Error Forecasts}

Theorem \ref{teo:csa_for} computes the minimum mean square error forecasts for a series generated by cross-sectional aggregation.

\begin{teo}\label{teo:csa_for}
	Let $x_t\sim CSA(a,b)$ defined as in (\ref{eq:csa_def}), and let $h\in\mathbb{N}$, then the minimum mean square error forecast $h$ periods ahead, $\hat{x}_{t+h}$, can be computed as\footnote{The Theorem assumes a Type II process analogous to samples of fractionally differenced processes, see \cite{Davidson2009}. In particular, it assumes $\nu_{j}=0,\ \ \forall j <0$.}
	\begin{equation}\label{eq:csa_for}
	\hat{x}_{T+h} = \sum_{j=h}^{T}\phi_{j}\nu_{T-j+h},
	\end{equation}
	where $\nu_{i}$ is computed as 
	\begin{equation*}\label{eq:csa_mat_nu}
	\nu_{i} = x_i - \sum_{j=1}^{i}{\phi_j \nu_{i-j}} \ \ \ \forall i\in{0,1,\cdots,T},
	\end{equation*}
	with $\phi_j$ as in Theorem \ref{teo:csa_gen}.
\end{teo}
\noindent Proof: See Appendix \ref{app:proofs}.

The theorem relies on the $MA(\infty)$ representation to construct the forecasts. Algorithms for computing the minimum mean square forecasts, (\ref{eq:csa_for}), are presented in Appendix \ref{app:codes}, and are available on the author's github repository at \href{https://github.com/everval}{github.com/everval/Nonfractional}.

Theorem \ref{teo:csa_for} allows for the construction of forecasts for the correct theoretical specification. Note that in the theorem we have assumed that we know the parameters of the the Beta distribution for the cross-sectionally aggregated process; for empirical applications we can use \pcite{Beran2010} method to estimate them. In the following, we will use these computations to assess the forecasting performance of high-order $AR$, and $ARFIMA$ models when working on long memory generated by cross-sectional aggregation.

\subsection{Forecasts With AR(p) Models}

This subsection computes the forecasting efficiency loss of an $AR(p)$ model when working with a $CSA(a,b)$ process. Theorem \ref{teo:arp} obtains the parameter estimates for an $AR(p)$ model fitted to a $CSA(a,b)$ process, and computes the efficiency loss for one-step ahead forecasts.

\begin{teo}\label{teo:arp}
	Let $x_t\sim CSA(a,b)$ defined as in (\ref{eq:csa_def}), and estimate an $AR(p)$ given by $(1-\alpha_1 L - \alpha_2 L^2 -\cdots -\alpha_p L^p )x_t=u_t$, then 
	\begin{equation*}
	\begin{bmatrix}	\alpha_1 \\ \alpha_2 \\ \vdots \\ \alpha_P	\end{bmatrix} = \begin{bmatrix}	1 &\gamma_{CSA(a,b)}(1) &\cdots &\gamma_{CSA(a,b)}(p-1) \\ \gamma_{CSA(a,b)}(1) & 1 &\cdots &\gamma_{CSA(a,b)}(p-2)\\ \vdots &\vdots &\ddots &\vdots \\  \gamma_{CSA(a,b)}(p-1)&\gamma_{CSA(a,b)}(p-2) &\cdots &1\end{bmatrix}^{-1}
	\begin{bmatrix}	\gamma_{CSA(a,b)}(1) \\ \gamma_{CSA(a,b)}(2) \\ \vdots \\ \gamma_{CSA(a,b)}(p) \end{bmatrix},
	\end{equation*}
	with $\gamma_{CSA(a,b)}(\cdot)$ defined as in \ref{ACF:CSA}.
	
	Furthermore, the one-step ahead forecast error variance relative to the minimum square error forecasts, $\zeta_{AR(p)}$, is given by 
	$$\zeta_{AR(p)} = \left(\frac{B(a,b-1)}{B(a,b)}\right)\left[ \left(1 + \sum_{i=1}^{p}{\alpha_i^2}\right) + 2\sum_{i=1}^{p}{\gamma_{CSA(a,b)}(i)\left(-\alpha_i+\sum_{j=1}^{p-i}{\alpha_j\alpha_{j+i}}\right)}\right].$$
\end{teo}
\noindent Proof: See Appendix \ref{app:proofs}.

As an example, consider estimating an $AR(1)$ model to forecast a cross-sectionally aggregated process. From Theorem \ref*{teo:arp}, the autoregressive parameter is $\alpha_1 = \gamma_{CSA(a,b)}(1) =  B(a+1/2,b-1)/B(a,b-1)$, while the one-step ahead forecast error variance is $\zeta_{AR(1)}=(B(a,b-1)/B(a,b))(1-\gamma_{CSA(a,b)}^2(1))$, which shows that the efficiency loss is a nonlinear function of both parameters of the $CSA(a,b)$ process. 

\begin{figure}[h!]
	\centering
	\includegraphics[scale=1]{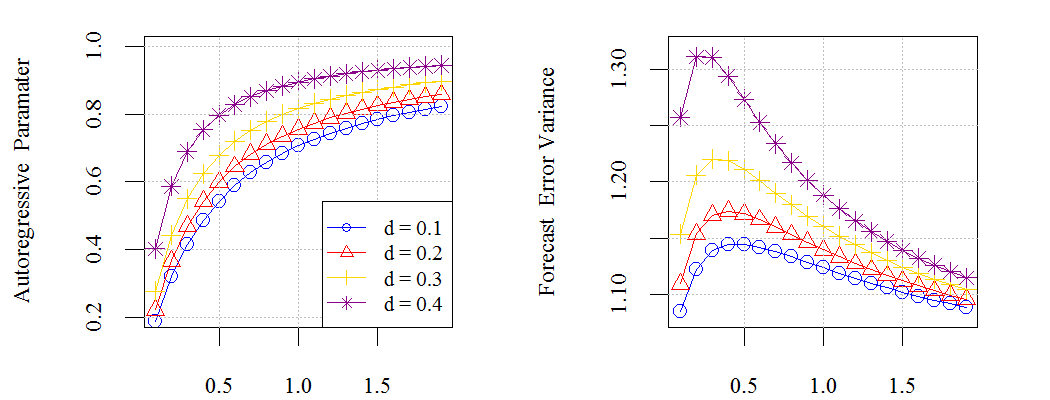}  
	\caption{Estimated autoregressive parameter and forecast error variance of fitted $AR(1)$ on $CSA(a,b)$ processes relative to the optimum forecasts while varying `$a$', for different degrees of memory.}
	\label{fig:For_AR1_CSA}
\end{figure}

To get a better sense of the nonlinearity, Figure \ref{fig:For_AR1_CSA} presents both the autoregressive parameter, $\alpha_1$, and the efficiency loss, $\zeta_{AR(p)}$, while varying the first parameter of the cross-sectional aggregated process, `$a$', with fixed $b=1.8$. On the one hand, the figure shows that as `$a$' increases, so does the autoregressive parameter, this is in line with the discussion above regarding the first parameter as a short memory regulator. On the other hand, the one-step ahead forecast error variance shows a maximum around $a=0.5$, with a forecast error variance of almost 15\%. That is, more than double the reported by \cite{Man2003} for an $I(0.1)$ process. Thus, the figures suggest that we could expect worse performance using an $AR(1)$ process to forecast a $CSA(a,b)$ process than a pure $I(d)$ one.

\begin{table}[ht!]
	\begin{footnotesize}
		\setlength{\tabcolsep}{3pt}
		\renewcommand{\arraystretch}{1.2} 
		\caption{One-step ahead forecast error variance of fitted $AR(1)$ and $AR(20)$ models on $CSA(a,b)$ processes relative to the optimum forecasts shown in (\ref{eq:csa_for}). The brackets show the associated long memory parameter $d$.} 
		\begin{center}
			\begin{tabular}{l|ccccc|ccccc}
				&\multicolumn{5}{c|}{$AR(1)$}&\multicolumn{5}{c}{$AR(20)$}\\
				\hline
				$a\backslash b\ [d]$ &1.8 [0.1] &1.6 [0.2] &1.4 [0.3]&1.2 [0.4]&1.1 [0.45]&1.8 [0.1] &1.6 [0.2] &1.4 [0.3]&1.2 [0.4]&1.1 [0.45]\\
				\hline		
				0.1	&		1.085		&		1.110		&		1.154		&		1.257		&		1.387		&	1.071	&	1.085	&	1.104	&	1.129	&	1.144	\\
				0.5	&		1.145		&		1.172		&		1.211		&		1.273		&		1.320		&	1.111	&	1.123	&	1.137	&	1.153	&	1.161	\\
				0.9	&		1.129		&		1.146		&		1.170		&		1.202		&		1.223		&	1.099	&	1.107	&	1.115	&	1.124	&	1.128	\\
				1.3	&		1.110		&		1.122		&		1.137		&		1.156		&		1.168		&	1.086	&	1.091	&	1.096	&	1.101	&	1.103	\\
				1.7	&		1.095		&		1.104		&		1.114		&		1.126		&		1.133		&	1.075	&	1.079	&	1.082	&	1.085	&	1.086	\\
				\hline
			\end{tabular}\label{tab:1fe_ar}
		\end{center}
	\end{footnotesize}
\end{table}

Table \ref{tab:1fe_ar} shows the one-step ahead forecast error variance of fitted $AR(1)$ models on $CSA(a,b)$ processes for different values of $a,b$. The losses are in line with the ones computed by \cite{Man2003} for the $I(d)$ case, increasing as the memory of the process increases. Moreover, as noted by \cite{Ray1993} and \cite{Man2003}, we can reduce the one-step ahead forecast error variance by allowing for more lags in the $AR$ specification. To get a sense of the improvements we can achieve, Table \ref{tab:1fe_ar} also presents the one-step ahead forecast error variance of $AR(20)$ models fitted to $CSA(a,b)$ processes for different values of $a,b$. As the table shows, increasing the order of the $AR$ model can greatly reduce the one-step ahead forecast error variance, particularly for larger degrees of long memory.

\subsection{Forecasts With Fractional Models}

This section studies the forecast performance of fractional models when working on cross-sectional aggregated processes. In particular, we compute the one-step ahead forecast error variance of the pure $I(d)$ model and of an $ARFIMA(1,d,0)$ allowing for more short term dynamics.

\begin{teo}\label{teo:fid}
	Let $x_t\sim CSA(a,b)$ defined as in (\ref{eq:csa_def}), then the one-step ahead forecast error variance of the $d=1-b/2$ fractional difference of the series, $\zeta_{I(d)}$, is given by 
	$$\zeta_{I(d)} = \gamma_z(0).$$
	
	Furthermore, estimate an $ARFIMA(1,d,0)$ by fitting an $AR(1)$ model to the $d=1-b/2$ fractional difference of the series, then the autoregressive parameter, $\alpha_{I}$, and the one-step ahead forecast error variance, $\zeta_{ARFIMA(1,d,0)}$, are given by 
	$$\alpha_{I} = \frac{\gamma_z(1)}{\gamma_z(0)},\ \ \ \ \zeta_{ARFIMA(1,d,0)} = \frac{\gamma_{z}(0)^2-\gamma_{z}(1)^2}{\gamma_{z}(0)^2}.$$
	
	Where in both expressions $\gamma_z(k)$ is given by
	$$\gamma_z(k) =  \frac{\gamma^*(k)}{B(a,b)}\left[B(a,b-1)\left(F_{1}(k)-1\right)+B(a+\frac{1}{2},q-1)F_{2}(k)\right],$$
	with 
	$$\gamma^*(k) = \sigma_{\varepsilon}^2\frac{\Gamma(1+2d)}{\Gamma(-d)\Gamma(1+d)}\frac{\Gamma(-d-k)}{\Gamma(1+d-k)},$$ and
	\begin{eqnarray*} 
		F_{1}(k) &:= &F\left[\left\{1,a,\frac{1-d+k}{2},\frac{-d+k}{2}\right\},\left\{a+b-1,\frac{2+d+k}{2},\frac{1+d+k}{2}\right\},1\right]+\\
		&&F\left[\left\{1,a,\frac{1-d-k}{2},\frac{-d-k}{2}\right\},\left\{a+b-1,\frac{2+d-k}{2},\frac{1+d-k}{2}\right\},1\right],\\
		F_{2}(k) &:= &\frac{-d+k}{1+d+k}*\\
		&&F\left[\left\{1,p+\frac{1}{2},\frac{1-d+k}{2},\frac{2-d+k}{2}\right\},\left\{a+b-\frac{1}{2},\frac{2+d+k}{2},\frac{3+d+k}{2}\right\},1\right]\\
		&&+\frac{-d-k}{1+d-k}*\\
		&&F\left[\left\{1,a+\frac{1}{2},\frac{1-d-k}{2},\frac{2-d-k}{2}\right\},\left\{a+b-\frac{1}{2},\frac{2+d-k}{2},\frac{3+d-k}{2}\right\},1\right],\\
	\end{eqnarray*}
	where $d=1-b/2$, and $F[\cdot]$ is the generalized hypergeometric function.
\end{teo}
\noindent Proof: See Appendix \ref{app:proofs}.

Table \ref{tab:1fe_arfima} presents the one-step ahead forecast error variance of fitted $ARFIMA(1,d,0)$, and $I(d)$ models on $CSA(a,b)$ processes for different values of $a,b$. The table also shows the estimated autoregressive parameter of the fitted $ARFIMA(1,d,0)$ model.

As the table shows, the pure fractional differenced process can be quite bad at forecasting $CSA(a,b)$ processes, specially for higher values of `$a$'. Once again, relating the first parameter functioning as a short memory regulator. Once we allow for more short term dynamics in the form of an $ARFIMA(1,d,0)$ model, the forecasting performance is much in line with the one from high order $AR$ models. Thus, it greatly helps to allow for some short term dynamics in the modelling scheme.

\begin{table}[h!]
	\begin{footnotesize}
		\setlength{\tabcolsep}{3pt}
		\renewcommand{\arraystretch}{1.2} 
		\caption{One-step ahead forecast error variance, and estimated autoregressive parameter in parenthesis, of fitted $I(d)$ and $ARFIMA(1,d,0)$ models on $CSA(a,b)$ processes. The brackets show the associated long memory parameter $d$.}
		\begin{center}
			\begin{tabular}{l|ccccc|ccccc}
				&\multicolumn{5}{c|}{$I(d)$}&\multicolumn{5}{c}{$ARFIMA(1,d,0)$}\\
				\hline
				$a\backslash b\ [d]$ &1.8 [0.1] &1.6 [0.2] &1.4 [0.3]&1.2 [0.4]&1.1 [0.45]&1.8 [0.1] &1.6 [0.2] &1.4 [0.3]&1.2 [0.4]&1.1 [0.45]\\
				\hline		
				0.1	&		1.077	&	1.084	&	1.112	&	1.158	&	1.186	&		1.072		&		1.083		&		1.103		&		1.131		&		1.147		\\	
				&			&		&		&		&		&	(	0.067	)	&	(	-0.019	)	&	(	-0.091	)	&	(	-0.153	)	&	(	-0.180	)	\\	
				0.5	&		1.345	&	1.253	&	1.202	&	1.176	&	1.169	&		1.127		&		1.132		&		1.138		&		1.147		&		1.152		\\	
				&			&		&		&		&		&	(	0.402	)	&	(	0.312	)	&	(	0.229	)	&	(	0.156	)	&	(	0.123	)	\\	
				0.9	&		1.615	&	1.435	&	1.318	&	1.240	&	1.212	&		1.118		&		1.121		&		1.123		&		1.125		&		1.125		\\	
				&			&		&		&		&		&	(	0.555	)	&	(	0.468	)	&	(	0.384	)	&	(	0.305	)	&	(	0.268	)	\\	
				1.3	&		1.880	&	1.611	&	1.431	&	1.309	&	1.263	&		1.104		&		1.107		&		1.108		&		1.107		&		1.106		\\	
				&			&		&		&		&		&	(	0.642	)	&	(	0.559	)	&	(	0.475	)	&	(	0.393	)	&	(	0.352	)	\\	
				1.7	&		2.138	&	1.778	&	1.538	&	1.374	&	1.312	&		1.093		&		1.096		&		1.097		&		1.095		&		1.093		\\	
				&			&		&		&		&		&	(	0.699	)	&	(	0.620	)	&	(	0.536	)	&	(	0.451	)	&	(	0.408	)	\\	
				\hline
			\end{tabular}\label{tab:1fe_arfima}
		\end{center}
	\end{footnotesize}
\end{table}

Comparing the results from Tables \ref{tab:1fe_ar} and \ref{tab:1fe_arfima}, we see that, as was the case for long memory series generated by the fractional difference operator, a high order $AR(p)$ can be a good model for forecasting long memory at short horizons. Yet, the bias-variance trade-off has to be assessed given the number of estimated parameters. As an alternative, the $ARFIMA(1,d,0)$ produces similar results while relying in only two parameters. 

\section{Conclusions}\label{Sec:Con}

Even though there is no theoretical argument linking the presence of long memory in the data with the fractional difference operator, fractionally differenced processes remain the most popular construction in the long memory time series literature. This may be due to the existence of efficient algorithms for their simulation and forecasting. Thus, this paper presents fast algorithms to generate long memory by a theoretically based mechanism that do not rely on the fractional difference operator, cross-sectional aggregation.

The cross-sectional aggregated process is then contrasted to the fractional difference operator. In particular, the paper analyses the antipersistent phenomenon. It is proven that, for negative degrees of memory, while the autocorrelations are negative by definition for the fractional difference operator, this restriction does not apply to the cross-sectional aggregated scheme. Furthermore, the paper shows that the lack of antipersistence for cross-sectional aggregated processes has implications for long memory estimators in the frequency domain which will be misspecified in general.

Moreover, this work evaluated the efficiency loss of using high-order $AR$ and $ARFIMA$ models to forecast long memory series generated by cross-sectional aggregation. It finds that, at short horizons, high-order $AR$ models beat a pure fractional difference process, $I(d)$, in terms of forecasting performance. Nonetheless, allowing for short term dynamic in the form of an $ARFIMA(1,d,0)$ model produces comparable forecast performance as high-order $AR$ models while relying on less parameters. 

The results of this paper can be used in the context of Monte Carlo simulations of long memory estimators and forecasts. Of particular interest is the analysis of long memory in inflation, one of the primer examples of cross-sectional aggregation producing long memory due to the way the data is computed. Moreover, the results allow for financial econometrics and climate econometrics models to incorporate long memory forecasts consistent with the theoretical argument posed by cross-sectional aggregation.

\section*{Acknowledgements}
I would like to thank Niels Haldrup for the careful reading of this article and all the insightful comments.


\section*{Appendix}
\renewcommand{\thesubsection}{\Alph{subsection}}

\subsection{Proofs}\label{app:proofs}

\subsubsection*{Proofs of Theorem \ref{teo:csa_gen}, and Theorem \ref{teo:ant}}

Lemma 1 of \cite{Haldrup2017}, which is in turn an application of the Central Limit Theorem, show that in the limit, the cross-sectionally aggregated process, Equation (\ref{eq:csa_def}), allows an $MA(\infty)$ representation. Using this representation, Theorem \ref*{teo:csa_gen} is an application of the circular convolution Theorem, see \cite{Jensen2014}.

Moreover, note that the spectral density for $x_t$ evaluated at the origin is given by 
$$f_X(0)=\frac{\sigma_\varepsilon^2}{2\pi} \left| \sum_{j=0}^\infty \phi_j \right|^2,$$
where $\phi_j$ are the coefficients of the $MA(\infty)$ representation given by $\phi_j=\frac{B(a+j,b)^{1/2}}{B(a,b)^{1/2}}$, with $B(\cdot,\cdot)$, the Beta function. Thus,
\begin{eqnarray*}
	f_X(0) &=&\frac{\sigma_\varepsilon^2}{2\pi} \left| \sum_{j=0}^\infty \phi_j \right|^2=\frac{\sigma_\varepsilon^2}{2\pi B(a,b)}  \left| \sum_{j=0}^\infty B(a+j,b)^{1/2} \right|^2=\frac{\sigma_\varepsilon^2}{2\pi B(a,b)} \left| \sum_{j=0}^\infty \frac{\Gamma(a+j)^{1/2}\Gamma(b)^{1/2}}{\Gamma(a+j+b)^{1/2}} \right|^2\\ &=&\frac{\sigma_\varepsilon^2\Gamma(b)}{2\pi B(a,b)} \left| \sum_{j=0}^\infty \frac{\Gamma(a+j)^{1/2}}{\Gamma(a+j+b)^{1/2}} \right|^2\approx \frac{\sigma_\varepsilon^2\Gamma(b)}{2\pi B(a,b)} \left| \sum_{j=0}^\infty j^{-b/2} \right|^2=\frac{\sigma_\varepsilon^2\Gamma(b)}{2\pi B(a,b)} \zeta(-b/2)^2\\
	&=&c_{a,b}<\infty,
\end{eqnarray*}
where $\Gamma(\cdot),\zeta(\cdot)$ are the Gamma and Euler–Riemann zeta function, respectively. Moreover, note that all terms in the expression are positive and thus $c_{a,b}>0$.

\subsubsection*{Proof of Theorem \ref{teo:csa_for}}

Assuming a Type II process, the Lemma shows that $x_t$ can be expressed as

\begin{equation}\label{eq:csa_mat}
\begin{bmatrix}	x_T \\ x_{T-1} \\ x_{T-2} \\ \vdots \\ x_0	\end{bmatrix} = \begin{bmatrix}	1 &\phi_1 &\phi_2 &\cdots &\phi_T \\ 0 & 1 &\phi_1 &\cdots &\phi_{T-1}\\ 0 & 0 & 1 &\cdots &\phi_{T-2}\\ \vdots &\vdots &\vdots &\ddots &\vdots \\ 0 & 0 &0 &\cdots &1\end{bmatrix} 
\begin{bmatrix}	\nu_T \\ \nu_{T-1} \\ \nu_{T-2} \\ \vdots \\ \nu_0 \end{bmatrix},
\end{equation}
with $\phi_t,\nu_t$ defined in Theorem \ref{teo:csa_gen}. Solving for $\nu_j$, this represents a system with an upper triangular matrix to which fast algorithms exist. Once we have solved it, the minimum mean square error forecast $h$ periods ahead can be easily obtained.

\subsubsection*{Proofs of Theorem \ref{teo:arp}, and Theorem \ref{teo:fid}} 

Let $x_t\sim CSA(a,b)$ and denote $\mathbf{\Phi} (\mathbf{z_{:,t}}):= \sum_{j=0}^\infty \phi_j z_{t-j}$ the function that applies the filter associated to the $MA(\infty)$ representation of the $CSA(a,b)$ process to $\mathbf{z_{:,t}}=(z_t, z_{t-1}, \cdots)$.  That is, 
$$x_t =  \sum_{j=0}^\infty \phi_j\nu_{t-j}:= \mathbf{\Phi}(\mathbf{\nu_{:,t}}),$$
where $\nu_j\sim N(0,\sigma^2_\varepsilon)$ and $\phi_j = \left(B(p+j,q)/B(p,q)\right)^{1/2},\ \forall j\in\mathbb{N}$. 

Let $u_t$ be the residual of a fitted $AR(p)$ process to $x_t$, i.e.,
$$u_t= (1-\alpha_1L-\alpha_2L^2-\cdots-\alpha_pL^p)x_t,$$
where $L$ is the lag operator and $\{\alpha_1,\cdots,\alpha_p\}$ are the fitted parameters. Hence,
\begin{eqnarray}\nonumber
u_t &= &(1-\alpha_1L-\alpha_2L^2-\cdots-\alpha_pL^p)x_t \\\nonumber
&= &(1-\alpha_1L-\alpha_2L^2-\cdots-\alpha_pL^p)\mathbf{\Phi}(\mathbf{\nu_{:,t}})\\\nonumber
&= &\mathbf{\Phi}(\mathbf{\nu_{:,t}})-\alpha_1\mathbf{\Phi}(\mathbf{\nu_{:,t-1}})-\alpha_2\mathbf{\Phi}(\mathbf{\nu_{:,t-2}})-\cdots-\alpha_p\mathbf{\Phi}(\mathbf{\nu_{:,t-p}})\\\nonumber
&=&\nu_t \phi_0+\nu_{t-1} (\phi_1-\alpha_1\phi_0)+\cdots+\nu_{t-s}(\phi_s-\cdots-\alpha_p\phi_{s-p})+\cdots\\\nonumber
&=&\sum_{j=0}^\infty\tilde{\phi}_j\nu_{t-j}, 
\end{eqnarray}
where 
$$\tilde{\phi_j}:=\left\{ \begin{array}{lr}
\phi_0 &j=0\\
\phi_1-\alpha_1\phi_0 &j=1\\
\phi_1-\alpha_1\phi_1-\alpha_2\phi_0  &j=2\\
\vdots&\vdots\\
\phi_j-\cdots-\alpha_p\phi_{j-p}&j\geq p
\end{array}\right.$$

Thus, 
\begin{eqnarray}\nonumber
var(u_t)\sigma_{\varepsilon}&=&\sum_{j=0}^\infty\tilde{\phi_j}^2\\\nonumber
&=&\phi_0^2 +(\phi_1-\alpha_1\phi_0)^2+\cdots+\sum_{j=p}^{\infty}(\phi_j -\alpha_1\phi_{j-1}-\cdots-\alpha_p \phi_{j-p})^2\\\nonumber
&=&\left(\sum_{j=0}^\infty \phi_j^2\right)(1+\alpha_1^2+\cdots+\alpha_p^2)+2\left(\sum_{j=1}^\infty\phi_j\phi_{j-1}\right)(-\alpha_1+\alpha_1\alpha_2+\cdots+\alpha_{p-1}\alpha_p)\\\nonumber
&&+2\left(\sum_{j=2}^\infty\phi_{j}\phi_{j-2}\right)(-\alpha_2+\alpha_1\alpha_3+\cdots+\alpha_{p-2}\alpha_p)+2\left(\sum_{j=p}^{\infty}\phi_j\phi_{j-p}\right)(-\alpha_p)\\\nonumber
&=& \left(\frac{B(a,b-1)}{B(a,b)}\right)\left[(1+\alpha_1^2+\cdots+\alpha_p^2)+2\gamma_{CSA(a,b)}(1)(-\alpha_1+\alpha_1\alpha_2+\cdots+\alpha_{p-1}\alpha_p)\right.\\\nonumber
&&\left.+\cdots+2\gamma_{CSA(a,b)}(p)(-\alpha_p)\right].
\end{eqnarray}			

Which shows that minimizing with respect to $\{\alpha_1,\cdots,\alpha_p\}$ yields the Yule-Walker system of equations, the result follows.

Finally, for Theorem \ref{teo:fid}, the first part regarding the autocorrelation function of a fractionally differenced cross-sectional aggregated process is Theorem 3 in \cite{Haldrup2017}. For the second part, the proof follows the same steps as for the pure $AR(p)$ case above replacing the autocorrelation function of the pure $CSA(a,b)$ process for the one of the fractionally differenced one.

\subsection{Code Implementation}\label{app:codes}

\subsubsection*{Nonfractional Memory Generation}

\textbf{R code}

\begin{lstlisting}[language=R]
csadiff = function(x, p, d){
iT = length(x)
n = nextn(2*iT - 1, 2)
k = 0:(iT-1)
b = (beta(p+k,d)/beta(p,d))^(1/2)
dx = fft(fft(c(x, rep(0, n - iT))) * 
fft(c(b, rep(0, n - iT))), inverse = T) / n;
return(Re(dx[1:iT]))
}
\end{lstlisting}

\noindent \textbf{Matlab code}

\begin{lstlisting}[language=Matlab]
function [cx] = csa_diff(x,p,q)
T = size(x,1);
np2 = 2.^nextpow2(2*T-1);
coefs = [1,( beta(p+(1:T-1),q) ./ ...
beta(p,q) ).^(1/2)];
cx = ifft(fft(x, np2).*fft(coefs', np2)); 
cx = cx(1:T, :);
end
\end{lstlisting}

\subsubsection*{Nonfractional Memory Forecasts}
\textbf{R code}

\begin{lstlisting}[language=R]
get_err_csa = function(x,p,q){
iT = length(x)
k = 0:(iT-1)
b = (beta(p+k,q)/beta(p,q))^(1/2)
zen = c(1, rep(0,iT-1))
ma = pracma::Toeplitz(b,zen)
vjs = solve(ma,x)

return(vjs)
}

forecast_csa = function(x,p,q,h){
iT = length(x)
k = seq(1,(iT),1)
b = (beta(p+k,q)/beta(p,q))^(1/2)

err = get_err_csa(x,p,q)
errv = rev(err)
fx = rep(0,h)
for(ii in 1:h){
fx[ii] = crossprod(b[ii:iT],errv[1:(iT-ii+1)])
}
return(fx)
}
\end{lstlisting}

\noindent \textbf{Matlab code}
\begin{lstlisting}[language=Matlab]
function [vjs] = get_err_csa(x,p,q)
T = size(x,1);
coefs = [1,( beta(p+(1:T-1),q) ./ beta(p,q) ).^(1/2)];
zen = [1,zeros(1,T-1)];
ma = toeplitz(coefs,zen);

vjs = ma\x;
end

function [fx] = forecast_csa(x,p,q,h)	
T = size(x,1);
coefs = [1,( beta(p+(1:T-1),q) ./ beta(p,q) ).^(1/2)];
coefs = fliplr(coefs);

err = get_err_csa(x,p,q);
fx = zeros(h,1);
for ii = 1:h
fx(ii,1) = coefs(1,1:T-ii)*err(ii+1:T,1);
end
end
\end{lstlisting}

\end{document}